\documentclass[12pt]{article}
\usepackage{amssymb}	% Allows use of AMS's mathematical symbols
\usepackage{amsmath}    % Allows use of AMS's mathematical commands
\usepackage{epsfig}
\usepackage{multirow}
\usepackage{longtable}
\usepackage{algorithmic}

\newcommand{\qed}{\nobreak \ifvmode \relax \else \ifdim\lastskip<1.5em \hskip-\lastskip \hskip1.5em plus0em minus0.5em \fi \nobreak \vrule height0.75em width0.5em depth0.25em\fi} 
\def\0{\bf \0}
\def\A{{\bf A}}

\def\0{{\bf 0}}

\def\R{\mathbb{R}}

\def\T{{\bf T}}

\def\a{{\bf a}}
\def\b{{\bf b}}
\def\c{{\bf c}}

\def\x{{\bf x}}
\def\y{{\bf y}}

\def\Tr{{\rm T}}
\def\T{{\rm T}}

\newtheorem{problem}{Problem}

\usepackage{color}
\usepackage[normalem]{ulem}

\begin{document}
\title{A collection of cycling problems in linear programming}
\author{
Yaguang Yang\thanks{US NRC, Office of Research, 
11555 Rockville Pike, Rockville, 20852. 
Email: yaguang.yang@verizon.net.} 
}

\date{\today}

\maketitle    % This command generates the title.

\begin{abstract}
This paper provides a set of cycling problems in linear
programming. These problems should be useful for 
researchers to develop and test new simplex algorithms.
As a matter of the fact, this set of problems is used 
to test a recently proposed double pivot
simplex algorithm for linear programming.
\end{abstract}


\newpage
%
%{\bf MSC classification:} 90C05 90C49.
% \newpage
% 

\section{Introduction}

Since Dantzig formulated linear programming (LP)
problem \cite{dantzig49} and invented the simplex 
method, LP has been extensively studied for more 
then seven decades.

There are two important classes of methods that 
have been developed for solving LP. The first one
is the simplex method which searches optimizer along
edges of the polyhedra from one vertex to the
next vertex. Many different pivot rules, such as 
Dantzig's most negative rule \cite{dantzig51}, the best
improvement rule \cite{jeroseow73}, Bland's least 
index pivoting rule \cite{bland77}, the steepest edge 
simplex rule \cite{fg92}, Zadeh's rule \cite{zadeh80}, 
among others \cite{tz93}, have been developed. 
The second class is the interior-point method 
which searches optimizer from interior of the 
polyhedra along the central path using either
line search technique \cite{wright97}
or arc-search technique \cite{yang20a}.
When solving LP problems, these two methods 
face different challenges. For the simplex method,
the cycling, due to degeneracy of several basic feasible
solutions for which the update cycles in a loop and stays
away from the optimal solution, can be a problem.
Although the degenerate/cycling 
case would seem (on mathematical grounds) to be
rare, this is not so in practice \cite{gale07}.
As a matter of fact, many benchmark problems
in Netlib repository have the degeneracy/cycling issue.
Therefore, several anti-cycling methods were 
developed. The degeneracy/cycling, however, is not a
severe problem for the interior-point method because
the latter does not search among vertices and its iterates in 
a carefully designed algorithm do not approach the
boundary before they approach to the optimal 
solution. A challenge for interior-point method is
that many LP problems do not have an interior
point, but this was resolved by using infeasible 
interior-point algorithms. 

An obvious merit of the interior-point method is that
many algorithms in this class are polynomial while 
none of the existing algorithms in simplex method has the property. 
However, the computational efficiency of the simplex
method has been demonstrated by many years of 
experience. Moreover, {\bf if} one finds a polynomial 
simplex algorithm, it will be strongly polynomial 
\cite{smale99}, a very attractive property that
interior-point algorithms do not have. To 
develop and test potentially better simplex algorithms, 
such as recently proposed double pivot algorithms 
\cite{ve18,yang20a}, it is desired to test them
on large size benchmark problems in Netlib repository.
Since many problems in Netlib repository have
the cycling issue, one may want to start working
on small size cycling LP problems before testing more 
challenging large size problems in Netlib repository. 
This motivated us to collect small size
cycling LP problems.\footnote{We realize that there 
exist several systematic methods to create simple 
cycling LP problems \cite{hm04,zornig08},
but significant efforts are required and we try to 
avoid repeating the efforts.}

Cycling is a phenomenon that the iterates move in 
a cycle, i.e., when a basic feasible solution is degenerate, 
after a few iterations using the simplex algorithm, it
may returns to a previously constructed tableau
\cite{zornig21} without improvement 
in objective function and they may stay away from the
optimal solution. Cycling happens when several
conditions are met. First, the problem has degenerate
basic feasible solutions\footnote{Having some degenerate 
basic feasible solutions is necessary for a linear programming 
problem to cycle in iterations, but it is not sufficient 
for cycling to occur.};  the initial basic feasible solution
is also a fact to determine if the cycling will happen;
moreover, cycling for one pivot rule does not
mean that it will happen for other pivot rules;
finally, when entering or leaving variables have a 
tie, different ways to break the tie also affect the 
cycling occurrence. Problems listed in this paper have a cycling
behavior for either Dantzig's pivot rule \cite{dantzig51}
or steepest edge rule \cite{gr77}. We will also
indicate which method is used to break the ties
for cycling to happen.

The remainder of the paper is organized as follows.
Section 2 provides the cycling problems for Dantzig's pivot
rule. Section 3 is the collection of the cycling problems
for steepest edge rule. Section 4 discusses the 
test result for a recently developed double pivot 
simplex algorithm against the cycling problems 
presented in this paper. Conclusions are summarized
in the last section.

\section{Cycling problems for Dantzig's pivot rule}

We divide problems in this group by either using the
least index rule or the largest pivot rule to break
ties in leaving variables.

\subsection{Cycling problems for the least index rule}

\begin{problem}
The first cycling example was given in \cite{gp07,gv04,hoffman53}
by Hoffman as follows\footnote{The original paper 
\cite{hoffman53} was not published. A detailed description of this
problem was given in \cite{gp07}. The form presented here was 
given in \cite{gv04}.}:
\[
\begin{array}{rl}
\min & -2.2361x_4 + 2x_5 + 4x_7 + 3.6180x_8 
+ 3.236x_9 + 3.6180x_{10} + 0.764x_{11} \\
s.t. & x_1  =1  \\
&  x_2+ 0.3090x_4 - 0.6180x_5 - 0.8090x_6 - 0.3820x_7 \\
&  + 0.8090x_8 + 0.3820x_9 + 0.3090x_{10}+0.6180x_{11} =0 \\
&  x3 + 1.4635x_4 + 0.3090x_5 + 1.4635x_6 - 0.8090x_7 \\
&  - 0.9045x_8 - 0.8090x_9 + 0.4635x_{10}+0.309x_{11} = 0 \\
& x_j \ge 0, \hspace{0.1in} j=1,\ldots, 11.
\end{array}
\]
The optimal solution is $x_1 = 1$ and $x_j =0$ ($j = 2 \ldots 11$) with
$obj=0$.
\end{problem}

\begin{problem}
The second cycling example was given in \cite{beale55}
by Beale as follows:
\[
\begin{array}{rl}
\min &  -3/4x_1 + 150x_2 - 1/50x_3 + 6x_4 \\
s.t. & 1/4 x_1 - 60x_2 - 1/25 x_3 + 9x_4 + x_5 = 0  \\
&  1/2 x_1 - 90x_2 - 1/50 x_3 + 3x_4 + x_6 = 0 \\
& x_3 + x_7 = 1 \\
& x_j \ge 0, \hspace{0.1in} j=1,\ldots, 7.
\end{array}
\]
Starting with initial base $[x_5,~x_6,~x_7]$ shows 
the cycling for Dantzig's pivot rule.
The optimal solution is $x_1 = 1/25$, $x_3 =1$, $x_5=3/100$ with
$obj=-1/20$.
\end{problem}

\begin{problem}
The following linear program problem was
given in \cite{chvatal83,gv04}\footnote{The following
three problems were given in \cite{gv04} without
providing initial points.}:
\[
\begin{array}{rl}
\min & 10 x_1-57x_2-9x_3-24x_4 \\
s.t. & 0.5x_1-5.5x_2-2.5x_3+9x_4 +x_5 =0  \\
& 0.5x_1 -1.5x_2 -0.5x_3 +x_4 + x_6 =0 \\
&  x_1 +x_7 =1 \\
& x_j \ge 0, \hspace{0.1in} j=1,\ldots, 7.
\end{array}
\]
In \cite{gv04}, the optimal solution was given as: $x_1=1$, 
$x_3=1$, $x_5=2$; 
and $\min =1$. As a matter of fact, au unbounded solutions 
of this problem is given as: $x_1=x_3=x_5=0$,
$x_7=1$, $x_2=\frac{9}{5.5} x_4$, and 
$x_6=\left( \frac{1.5*9}{5.5}-1\right)x_4$; while
$x_4 \rightarrow \infty$, the objective goes to
$-\infty$. Another unbounded solution is given as:
$x_1=x_5=x_7=0$, $x_2=c$, $x_3=0.5c$, $x_4=0.75c$,
and $x_6=c$ with $c \ge 0$; as $c \rightarrow \infty$,
the objective goes to $-\infty$.
\end{problem}

\begin{problem}
The following linear program problem was presented 
in \cite{yg65}: 
\[
\begin{array}{rl}
\max & x_3- x_4 +x_5-x_6 \\
s.t. & x_1 +x_3 -2x_4 -3x_5 +4x_6=0  \\
& x_2 +4x_3 -3x_4 -2x_5+ x_6 =0 \\
&  x_3+x_4+x_5+x_6 +x_7 =1 \\
& x_j \ge 0, \hspace{0.1in} j=1,\ldots, 7.
\end{array}
\]
The optimal solution was given in \cite{gv04} as: $x_1=3$, $x_2=2$, 
$x_5=1$; and $\max =1$. For this problem, there is another 
optimal solution which is given as $x_1=\frac{5}{3}$, 
$x_3=\frac{1}{3}$, and $x_5=\frac{2}{3}$ with $obj =1$.
\end{problem}

\begin{problem}
The following linear program problem was presented 
in \cite{yg65}:
\[
\begin{array}{rl}
\min &  -x_3 + x_4 - x_5 + x_6  \\ 
s.t. & 
    x_1 + 2x_3 - 3x_4 - 5x_5 + 6x_6 = 0 \\
& x_2 + 6x_3 - 5x_4 - 3x_5 + 2x_6 = 0  \\
&    3x_3 + x_4 + 2x_5 + 4x_6 + x_7 = 1 \\
& x_j \ge 0, \hspace{0.1in} j=1,\ldots, 7.
\end{array}
\]
The optimal solution is given as: $x_1=2.5$, $x_2=1.5$, 
$x_5=0.5$; and $obj =-0.5$. 
\end{problem}

\begin{problem}
The following linear program problem was presented
in \cite{bt97}:
\[
\begin{array}{rl}
\min &  - 2x_4 - 3x_5 + x_6 + 12x_7  \\ 
s.t. & 
   x_1 - 2x_4 - 9x_5 + x_6 + 9x_7 = 0 \\
& x_2 + 1/3x_4 + x_5 - 1/3x_6 - 2x_7 = 0  \\
&   x_3 + 2x_4 + 3x_5 - x_6 - 12x_7 = 2 \\
& x_j \ge 0, \hspace{0.1in} j=1,\ldots, 7.
\end{array}
\]
Starting with initial base $[x_1,~x_2,~x_3]$ shows 
the cycling for Dantzig's pivot rule.
The optimal solution is given as: $x_1=2$, $x_4=2$, 
$x_6=2$; and $obj=-2$. 
\end{problem}

\begin{problem}
The following linear program problem was presented
in \cite{ms69}:
\[
\begin{array}{rl}
\min &   - 0.4x_5 - 0.4x_6 + 1.8x_7  \\ 
s.t. & 
   x_1 + 0.6x_5 - 6.4x_6 + 4.8x_7 = 0   \\
& x_2 + 0.2x_5 - 1.8x_6 + 0.6x_7 = 0  \\
& x_3 + 0.4x_5 - 1.6x_6 + 0.2x_7 = 0  \\
& x_4 + x_6 = 1 \\
& x_j \ge 0, \hspace{0.1in} j=1,\ldots, 7.
\end{array}
\]
Starting with initial base $[x_1,~x_2,~x_3,~x_4]$ 
shows the cycling for Dantzig's pivot rule.
The optimal solution is given as: $x_1 = 4$, $x_2 = 1$, 
$x_5 = 4$, $x_6 = 1$ with $obj= -2$. 
\end{problem}

\begin{problem}
The following linear program problem was presented
in \cite{gv04,solow84}\footnote{The following
three problems were given in \cite{gv04} without
providing initial points}:
\[
\begin{array}{rl}
\min &   - 2x_3 -2x_4 + 8x_5 + 2x_6  \\ 
s.t. & 
   x_1 - 7x_3 - 3x_4 + 7x_5 + 2x_6 = 0   \\
& x_2 + 2x_3 + x_4 - 3x_5 - x_6 = 0 \\
& x_j \ge 0, \hspace{0.1in} j=1,\ldots, 6.
\end{array}
\]
The optimal solution is given as: for $c \ge 0$, $x_1=x_4=x_6=c$,
and the rest variables are zeros 
(while a special case is given in \cite{gv04} that all variables are zeros) 
with $obj=0$. 
\end{problem}

\begin{problem}
The following linear program problem was presented
in \cite{sierksma96}:
\[
\begin{array}{rl}
\min &   -3x_1 + 80x_2 - 2x_3 + 24x_4  \\ 
s.t. & 
x_1 - 32x_2 - 4x_3 + 36x_4 + x_5 = 0   \\
& x_1 - 24x_2 - x_3 + 6x_4 + x_6 = 0 \\
& x_j \ge 0, \hspace{0.1in} j=1,\ldots, 6.
\end{array}
\]
The optimal solution is given as: $x_1=1.8c$, $x_2=x_6=0$,
$x_3=3c$, $x_4=0.2c$, $x_5=3c$, as $c \rightarrow \infty$,
the $obj \rightarrow -\infty$, unbounded solution. 
\end{problem}

\begin{problem}
The following linear program problem was presented
in \cite{nt93}:
\[
\begin{array}{rl}
\min &   -3x_2 + x_3- 6x_4 - 4x_6  \\ 
s.t. & x_1 + x_2 + 1/3x_5 + 1/3x_6 = 2   \\
&  9x_2 + x_3 - 9x_4 - 2x_5 - 1/3x_6 + x_7 = 0 \\
&  x_2 + 1/3x_3 - 2x_4 - 1/3x_5 - 1/3x_6 + x_8 = 2 \\
& x_j \ge 0, \hspace{0.1in} j=1,\ldots, 8.
\end{array}
\]
The optimal solution is given as: $x_1=x_2=x_3=x_5=0$,
$x_4=c/9$, $x_6=6$, $x_7=c$, $x_8=2c/9$, as $c \rightarrow \infty$,
the $obj \rightarrow -\infty$, unbounded solution. 
\end{problem}

\begin{problem}
The following linear program problem was presented
in \cite{zornig08}:
\[
\begin{array}{rl}
\min &    -3x_1 -59/20x_2 +50x_3 +2/5x_4  \\ 
s.t. &  1/40x_1 +1/400x_2 + 3x_3+    2x_4    +x_5 =0 \\
&   1/20x_1 +9/200x_2 -1/2x_3 +2/25x_4  +x_6 =0  \\
& x_j \ge 0, \hspace{0.1in} j=1,\ldots, 6.
\end{array}
\]
Starting with initial base $[x_5,~x_6]$ and
breaking the tie by using the first pivot,
the problem is cycling for Dantzig's rule.
The optimal solution is given as: all variables
are zero with $obj =0$. 
\end{problem}

\begin{problem}
The following linear program problem was presented
in \cite{zornig08}:
\[
\begin{array}{rl}
\min &   -14x_1 +25x_2 -7/20x_3 +20x_4  \\ 
s.t. &  x_1 -2x_2 -1/10x_3 +5x_4    +x_5 =0 \\
&   7/10x_1 -3/10x_2-1/100x_3 +19/50x_4  +x_6 =0  \\
& x_j \ge 0, \hspace{0.1in} j=1,\ldots, 6.
\end{array}
\]
Starting with initial base $[x_5,~x_6]$ and
breaking the tie by using the first pivot,
the problem is cycling for Dantzig's rule.
The optimal solution is given as: $x_2=x_4=x_6=0$,
$x_1=c/70$, $x_3=c$, $x_5=6c/70$, as $c \rightarrow \infty$,
the $obj \rightarrow - \infty$, unbounded solution.  
\end{problem}

\begin{problem}
The following linear program problem was presented
in \cite{zornig08}:
\[
\begin{array}{rl}
\min &   -1/2x_1 -2/5x_2 +5x_3 +1/5x_4  \\ 
s.t. &  1/40x_1 -1/100x_2 + 3x_3 +2x_4   +x_5 =0 \\
&   1/20x_1 +1/50x_2 +1/50x_3 +2/25x_4 +x_6 =0  \\
& x_j \ge 0, \hspace{0.1in} j=1,\ldots, 6.
\end{array}
\]
Starting with initial base $[x_5,~x_6]$ and
breaking the tie by using the first pivot,
the problem is cycling for Dantzig's rule.
The optimal solution is given as: all variables are
zeros with $obj=0$. 
\end{problem}

\begin{problem}
The following linear program problem was presented
in \cite{zornig08}:
\[
\begin{array}{rl}
\min &  -14x_1 +25x_2 -7/20x_3 +20x_4  \\ 
s.t. &  x_1 -2x_2 -1/10x_3 +5x_4   +x_5 =0 \\
&   7/10x_1 -3/10x_2 -1/100x_3 +19/50x_4 +x_6 =0  \\
&   x_1 +x_2 +x_3 +x_4 +x_7 =5  \\
&   x_1 +2x_2 +3x_3 +x_4 +x_8 =10  \\
& x_j \ge 0, \hspace{0.1in} j=1,\ldots, 8.
\end{array}
\]
Starting with initial base $[x_5,~x_6,x_7,~x_8]$ and
breaking the tie by using the first pivot,
the problem is cycling for Dantzig's rule.
The optimal solution is given as: $x_1=10/211$, $x_2=0$, 
$x_3=700/211$, $x_4=0$, $x_5=60/211$, $x_6=0$,
$x_7=345/211$, and $x_8=0$ with $obj=-385/211$. 
\end{problem}

\begin{problem}
The following linear program problem was presented
in \cite{zornig08}:
\[
\begin{array}{rl}
\min &  -1/100x_1+1/100x_2 -9/1000x_3 +3/200x_4 -1/500x_5 +3/20x_6  \\ 
s.t. &  1/20x_1 -100x_2 -2/5x_3 -100x_4 -x_5 +65x_6    +x_7 =0 \\
&   9/10x_1 -x_2 +3/5x_3 -3/2x_4 -1/100x_5 +1/100x_6  +x_8 =0  \\
& x_j \ge 0, \hspace{0.1in} j=1,\ldots, 8.
\end{array}
\]
Starting with initial base $[x_7,~x_8]$ and
breaking the tie by using the first pivot,
the problem is cycling for Dantzig's rule.
The optimal solution is given as: $x_1=x_2=x_4=x_6=x_8=0$, 
$x_3=c/60$, $x_5=149c/150$, $x_7=c$, as $c \rightarrow \infty$,
$obj  \rightarrow - \infty$, unbounded. 
\end{problem}

\begin{problem}
The following linear program problem was presented
in \cite{zornig08}:
\[
\begin{array}{rl}
\min &  -3/100x_1 -1/100x_2 -1/100x_3 +x_4 +3/10x_5 +1/10x_6 +2/125x_7 +1/2x_8
  \\ 
s.t. &  1/10x_1 -100x_2 -13x_3 -3/20x_4 -6x_5 +23/100x_6 +1/100x_7 +10x_8   +x_9 =0 \\
&   1/2x_1 +3/5x_2 +2/25x_3 -8x_4 -5x_5 -13/10x_6 -2/5x_7 +1/10x_8  +x_{10} =0  \\
& x_j \ge 0, \hspace{0.1in} j=1,\ldots, 10.
\end{array}
\]
Starting with initial base $[x_9,~x_{10}]$ and
breaking the tie by using the first pivot,
the problem is cycling for Dantzig's rule.
The optimal solution is given as:  $x_2=x_4=x_5=x_6=x_8=x_{10}=0$, 
$x_1=0.786601106330670c$, $x_3=0.083743085433313c$, $x_7=x_9=c$, 
as $c \rightarrow \infty$, $obj  \rightarrow - \infty$, unbounded. 
\end{problem}

\begin{problem}
The following linear program problem was presented
in \cite{zornig08}:
\[
\begin{array}{rl}
\min &   -14x_1 +25x_2 -7/20x_3 +20x_4  \\ 
s.t. &  x_1 -2x_2 -1/10x_3 +5x_4    +x_5 =0 \\
&   7/10x_1 -3/10x_2-1/100x_3 +19/50x_4  +x_6 =0  \\
& x_j \ge 0, \hspace{0.1in} j=1,\ldots, 6.
\end{array}
\]
Starting with initial base $[x_5,~x_6]$ and
breaking the tie by using the first pivot,
the problem is cycling for steepest edge rule.
The optimal solution is given as: $x_2=x_4=x_6=0$, 
$x_1=c$, $x_3=70c$, $x_5=6c$, 
as $c \rightarrow \infty$, $obj  \rightarrow - \infty$, unbounded.  
\end{problem}

\begin{problem}
The following linear program problem was presented
in \cite{zornig08}:
\[
\begin{array}{rl}
\min &   -1/4x_4 -23/100x_5 -4849/20000x_6 +21/50x_7 \\
&			-123/2000x_8 -1809/10^6x_9 +4511/1250x_{10}  \\ 
s.t. &  x_1 +2x_4 +6/5x_5 +13/10x_6 +1/100x_7 +7/10x_8 +1/1000x_9 +3/50x_{10} = 0  \\
&   x_2 +7/5x_4 +13/10x_5 +34/25x_6 +1/20x_7 +6/5x_8 +13/10000x_9 +23/10x_{10} = 0   \\
&   x_3 -4x_4 -3/2x_5 -17/10x_6 -28/5x_7 -2x_8 -1/100x_9 +15x_{10} = 0  \\
& x_j \ge 0, \hspace{0.1in} j=1,\ldots, 10.
\end{array}
\]
Starting with initial base $[x_1,~x_2,~x_3]$ and
breaking the tie by using the first pivot,
the problem is cycling for Dantzig's rule.
The optimal solution is given as: all variables are
zeros  with $obj=0$. 
\end{problem}

\subsection{Cycling problems for the largest pivot}

\begin{problem}
The following linear program problem was presented
in \cite{ms69}:
\[
\begin{array}{rl}
\min &   -x_3 + 7x_4 + x_5+ 2x_6  \\ 
s.t. &    x_1 +0.5 x_3 - 5.5x_4 - 2.5x_5 + 9x_6 = 0   \\
&  x_2 + 0.5x_3 - 1.5x_4 - 0.5x_5 + x_6 = 0  \\
& x_j \ge 0, \hspace{0.1in} j=1,\ldots, 6.
\end{array}
\]
Starting with initial base $[x_1,~x_2]$ and
breaking the tie by using the largest pivot,
(as is normal for numerical stability),
the problem is cycling for Dantzig's pivot rule.
The optimal solution is given as: for $c \ge 0$, $x_1=2c$, $x_3=x_5=c$,
and the rest variables are zeros (while \cite{ms69} gives a 
special solution that all variables are zero) with $obj =0$.
\end{problem}

\begin{problem}
The following linear program problem was also presented
in \cite{hm04} which introduces two additional constraints
to make it a bounded problem:
\[
\begin{array}{rl}
\min &   -2.3x_1 - 2.15x_2 + 13.55x_3 + 0.4x_4  \\ 
s.t. &  0.4x_1 + 0.2x_2 - 1.4x_3 - 0.2x_4 +x_5 = 0  \\
&  -7.8x_1 - 1.4x_2 + 7.8x_3 + 0.4x_4 + x_6 = 0 \\
&  x_1 +  x_7 = 1  \\
&  x_2 +  x_8 = 1  \\
& x_j \ge 0, \hspace{0.1in} j=1,\ldots, 8.
\end{array}
\]
Starting with initial base $[x_5,~x_6,~x_7,~x_8]$ and
breaking the tie by using the largest pivot,
(as is normal for numerical stability),
the problem is cycling for Dantzig's pivot rule.
The optimal solution is given as: $x_1 = 1$, $x_2 = 1$, 
$x_4 = 3$, $x_6 = 8$ , and $x_3=x_5=x_7=x_8=0$
with $obj= -3.25$.
\label{p14}
\end{problem}

\begin{problem}
The following linear program problem was presented
in \cite{hm04}:
\[
\begin{array}{rl}
\min &   -2.3x_1 - 2.15x_2 + 13.55x_3 + 0.4x_4  \\ 
s.t. &  0.4x_1 + 0.2x_2 - 1.4x_3 - 0.2x_4 +x_5 = 0  \\
&  -7.8x_1 - 1.4x_2 + 7.8x_3 + 0.4x_4 + x_6 = 0 \\
& x_j \ge 0, \hspace{0.1in} j=1,\ldots, 6.
\end{array}
\]
Starting with initial base $[x_5,~x_6]$ and
breaking the tie by using the largest pivot,
(as is normal for numerical stability),
the problem is cycling for Dantzig's pivot rule.
The optimal solution is given as: $x_1=x_3=x_5=0$, 
$x_2=c$, $x_4=c$, $x_6=c$, 
as $c \rightarrow \infty$, $obj  \rightarrow - \infty$, unbounded. 
\end{problem}

\begin{problem}
The following linear program problem was presented
in \cite{zornig08}:
\[
\begin{array}{rl}
\min &  -1/100x_1 -1/1000x_2 -1/2000x_3 +63/100x_4  \\ 
s.t. &  20x_1 +7/100x_2 -7/100x_3 + 100x_4  +x_5 = 0  \\
& -100x_1 -3/10x_2 -1/100x_3 +1/4x_4  + x_6 = 0 \\
& x_j \ge 0, \hspace{0.1in} j=1,\ldots, 6.
\end{array}
\]
Starting with initial base $[x_5,~x_6]$ and
breaking the tie by using the largest pivot,
(as is normal for numerical stability),
the problem is cycling for Dantzig's pivot rule.
The optimal solution is given as: $x_1=x_4=x_5=0$, 
$x_2=c$, $x_3=c$, $x_6=0.31c$, 
as $c \rightarrow \infty$, $obj  \rightarrow - \infty$, unbounded. 
\end{problem}

\begin{problem}
The following linear program problem was presented
in \cite{zornig08}:
\[
\begin{array}{rl}
\min &  -1/250x_1 +1/200x_2 -1/1000x_3 +3/20x_4+1/2000x_5 +1/10x_6  \\ 
s.t. &  x_1 -100x_2 -1/25x_3 -1/250x_4 -1/25x_5 +3/2x_6 +x_7 = 0  \\
& 1/2x_1 -1/10x_2 +1/1000x_3 -1/2x_4 -2/25x_5 +1/2x_6  + x_8 = 0 \\
& x_j \ge 0, \hspace{0.1in} j=1,\ldots, 8.
\end{array}
\]
Starting with initial base $[x_7,~x_8]$ and
breaking the tie by using the largest pivot,
(as is normal for numerical stability),
the problem is cycling for Dantzig's pivot rule.
The optimal solution is given as: $x_2=x_4=x_6=x_7=x_8=0$, 
$x_1=0.054c$, $x_3=c$, $x_5= 0.35c$, 
as $c \rightarrow \infty$, $obj  \rightarrow - \infty$, unbounded. 
\end{problem}

\begin{problem}
The following linear program problem was presented
in \cite{zornig08}:
\[
\begin{array}{rl}
\min & -8x_1+25x_2 -39/5x_3 +3182/5x_4 -37/25x_5 \\
 &     +4713/1000x_6 -2447/2500x_7 +247367/5000x_8  \\ 
s.t. &  11/25x_1 -50x_2 -3x_3 +12x_4 -8x_5 +1/2x_6 -x_7 +50x8 +x_9 = 0  \\
& 2/5x_1 +x_2 +9/100x_3 -5/2x_4 -1/100x_5 -1/50x_6 -1/2500x_7 +1/100x_8  + x_{10} = 0 \\
& x_j \ge 0, \hspace{0.1in} j=1,\ldots, 10.
\end{array}
\]
Starting with initial base $[x_9,~x_{10}]$ and
breaking the tie by using the largest pivot,
(as is normal for numerical stability),
the problem is cycling for Dantzig's pivot rule.
The optimal solution is given as: $x_1=x_2=x_4=x_5=x_6=x_7=x_8=x_{10}=0$, 
$x_3=0.004385964912281c$, $x_7=0.986842105263158c$, $x_9= c$, 
as $c \rightarrow \infty$, $obj  \rightarrow - \infty$, unbounded. 
\end{problem}

\begin{problem}
The following linear program problem was presented
in \cite{zornig08}:
\[
\begin{array}{rl}
\min & -x_1 +340x_2 +71/5x_3 +107/2x_4 -7/40x_5 +1469/500x_6 -99/100x_7 +15627/100x_8  \\ 
s.t. &  2/5x_1 -145x_2 -39/5x_3 +6/5x_4 -11/50x_5 +11/5x_6 -11/50x_7 -1/2x_8 +x_{9} = 0  \\
& 9/25x_1 +28x_2 +10x_3 -120x_4 -1/2x_5 -2x_6 -x_7 +100x_8  + x_{10} = 0 \\
& x_j \ge 0, \hspace{0.1in} j=1,\ldots, 10.
\end{array}
\]
Starting with initial base $[x_9,~x_{10}]$ and
breaking the tie by using the largest pivot,
(as is normal for numerical stability),
the optimal problem is cycling for Dantzig's pivot rule.
The solution is given as: $x_2=x_4=x_6=x_8=x_{9}=0$, 
$x_1=c$, $x_3=0.003483870967742c$, $x_5=0.789677419354839c$, 
$x_7=x_{10}= c$, 
as $c \rightarrow \infty$, $obj  \rightarrow - \infty$, unbounded. 
\end{problem}

\begin{problem}
The following linear program problem was presented
in \cite{zornig08}:
\[
\begin{array}{rl}
\min & -1/10x_1 -1/20x_2 -2/25x_3 +1/10x_4 -1/25x_5 -2/25x_6 +33/200x_7  \\ 
s.t. &  2/5x_1 -13/10x_2 +5x_3 +8/5x_4 +1/5x_5 +3/5x_6 +3/2x_7 +x_{8} = 0  \\
& 12/5x_1 +2/5x_2 +2x_3 -13/20x_4 +1/25x_5 +7/10x_6 -1/2x_7 + x_{9} = 0 \\
& -25x_1 -8/5x_2 -24x_3 -1/5x_4 -2x_5 -3x_6 +5x_7 +x_{10} = 0   \\
& x_j \ge 0, \hspace{0.1in} j=1,\ldots, 10.
\end{array}
\]
Starting with initial base $[x_8,~x_9,~x_{10}]$ and
breaking the tie by using the largest pivot,
(as is normal for numerical stability),
the problem is cycling for Dantzig's pivot rule.
The optimal solution is given as: $x_1=x_3=x_6=x_7=x_8=x_{9}=0$, 
$x_2=0.259775040171398c$, $x_4=0.176754151044456c$, $x_5=0.274504552758436c$, $x_{10}= c$, 
as $c \rightarrow \infty$, $obj  \rightarrow - \infty$, unbounded. 
\end{problem}

\section{Cycling problems for steepest edge rule}

The following problems shows cycling behavior for steepest edge rule.

\begin{problem}
This linear program problem was presented in 
\cite{hm04} which introduces two additional constraints
in Problem \ref{p14} to make it a bounded problem:
\[
\begin{array}{rl}
\min &   -1.0x_1 - 1.75x_2 + 12.25x_3 + 0.5x_4  \\ 
s.t. &  0.4x_1 + 0.2x_2 - 1.4x_3 - 0.2x_4 +x_5 = 0  \\
&  -7.8x_1 - 1.4x_2 + 7.8x_3 + 0.4x_4 + x_6 = 0 \\
&   - 20x_2 + 156x_3 + 8.0x_4 + x_7 = 1 \\
& x_j \ge 0, \hspace{0.1in} j=1,\ldots, 7.
\end{array}
\]
Starting with initial base $[x_5,~x_6,~x_7]$ and
breaking the tie by using the larger pivot,
(as is normal for numerical stability),
the problem is cycling for steepest edge rule.
The optimal solution is given as: $x_1=x_3=x_5=0.070512820513462c$, 
$x_2=x_4=x_6=x_7= c$, as $c \rightarrow \infty$, $obj  \rightarrow - \infty$, unbounded. 
\end{problem}

\begin{problem}
The following linear program problem was presented
in \cite{web02} which is verified by the author:
\[
\begin{array}{rl}
\min &   -10x_1 +57x_2 + 9x_3 + 24x_4  \\ 
s.t. &  0.5x_1 -5.5x_2   -2.5x_3+    9x_4    +x_5 =0 \\
  &     0.5x_1   -1.5x_2   -0.5x_3+ x_4       +x_6 =0  \\
  &     x_1    +x_2   +x_3    +x_4    +x_7  =1  \\
& x_j \ge 0, \hspace{0.1in} j=1,\ldots, 7.
\end{array}
\]
Starting with initial base $[x_5,~x_6,~x_7]$ and
breaking the tie by using the first pivot,
the problem is cycling for steepest edge rule.
The optimal solution is given as: $x_1 = 0.5$, $x_3 = 0.5$, 
$x_5 = 1$, and $x_2=x_4=x_6=x_7=0$ with $obj= -0.5$. 
\end{problem}

\begin{problem}
This linear program problem was also presented
in \cite{zornig08}:
\[
\begin{array}{rl}
\min &   -7/100x_1 -3/50x_2 +11/50x_3 +133/12500x_4  \\ 
s.t. &  33/25x_1 +1/2x_2 -53/25x_3 -9/20x_4 +x_5 = 0  \\
&  -44x_1 -24/5x_2 +6x_3 +33/100x_4 + x_6 = 0 \\
&  1319/1000x_1 -40x_2 +2123/10x_3 +1173/100x_4 + x_7 = 0 \\
& x_j \ge 0, \hspace{0.1in} j=1,\ldots, 7.
\end{array}
\]
Starting with initial base $[x_5,~x_6,~x_7]$ and
breaking the tie by using the larger pivot,
(as is normal for numerical stability),
the problem is cycling for steepest edge rule.
The optimal solution is given as: $x_1=x_3=0$, 
$x_2=0.264476614699332c$, $x_4= 0.816629547141797c$, 
$x_5=0.235244988864143c$, $x_6=x_7=c$, 
as $c \rightarrow \infty$, $obj  \rightarrow - \infty$, unbounded. 
\end{problem}

\begin{problem}
This linear program problem was also presented
in \cite{zornig08}:
\[
\begin{array}{rl}
\min &    36x_1 -3/5x_2 +20x_3 +1/4x_4 -1/20x_5 -1/20x_6  \\ 
s.t. &  2x_1 +1/5x_2 -5x_3 -9/10x_4 +x_5 +23/1000x_6 +x_7 = 0  \\
&  -41x_1 -6/5x_2 +12x_3 +1/5x_4 -14/5x_5 -1/500x_6  + x_8 = 0 \\
&  165000x_1 +2600x_2 +9600x_3 +125x_4 -100x_5 -300x_6 + x_9 = 0 \\
& x_j \ge 0, \hspace{0.1in} j=1,\ldots, 9.
\end{array}
\]
Starting with initial base $[x_7,~x_8,~x_9]$ and
breaking the tie by using the larger pivot,
(as is normal for numerical stability),
the problem is cycling for steepest edge rule.
The optimal solution is given as: $x_1=x_3=x_5=x_7=x_9=0$, 
$x_2=0.112949260042283c$, $x_4= 0.050655391120507c$, 
$x_6=c$, $x_8=0.127408033826638c$, 
as $c \rightarrow \infty$, $obj  \rightarrow - \infty$, unbounded. 
\end{problem}

\section{Preliminary test for a double pivot algorithm}

Although Bland's rule can prevent cycling from happening,
it is not as efficient as the Dantzig's pivot rule \cite{gale07}.
A double pivot simplex algorithm was recently proposed
\cite{yang20a} aiming at dealing with degeneracy/cycling
problem and improving the efficiency of Dantzig's pivot
rule. This algorithm updates two pivots at a time which 
is different from the traditional simplex algorithms.
An important feature of the algorithm is that one of the
pivot takes the longest step among all possible entering
variables. This feature makes it possible to avoid cycling 
problems. The promising result of the double pivot simplex 
algorithm for large size LP problems has been discussed
in \cite{yang20a}. 
We consider the primal linear programming 
problem in the standard form:
\begin{eqnarray}
\begin{array}{cl}
\min  &  \c^{\Tr}\x, \\ 
\mbox{\rm subject to} 
&  \A\x=\b, \hspace{0.1in} \x \ge \0,
\end{array}
\label{LP}
\end{eqnarray}
where $\A=[\a_1, \ldots, \a_n] \in {\R}^{m \times n}$, 
$\b \in {\R}^{m}$, 
$\c \in {\R}^{n}$ are given, and $\x \in {\R}^n$  is the 
vector to be optimized. To save space, we will write the 
column vector $\x =[\x_1^{\Tr}, \x_2^{\Tr}]^{\Tr}$ as 
$\x=(\x_1, \x_2)$. We denote by $\a_i$, the $i$th column of 
$\A$, for $i \in \mathcal{I} = \{ 1, \ldots,n \}$, where 
the subscript $i$ is the index of column of $\A$. We 
denote by $B \subset \mathcal{I}$ the index set with 
cardinality $| B | = m$ and $N = \mathcal{I} \setminus B$ 
the complementary set of $B$ with cardinality $| N |=n-m$ 
such that matrix $\A$ and vector $\x$ can be partitioned as 
$\A=[ \A_B, \A_N]$ and $\x =(\x_B, \x_N)$, moreover the
columns of $\A_B$ are linearly independent and
$\A_B \x_B = \b$, hence $\x_N=\0$. We call this
$\x=(\x_B,\0) \ge \0$ as the basic feasible solution. Similarly, we
have the partition $\c=(\c_B, \c_N)$. 

Using the $B-N$ partition, we can rewrite the problem (\ref{LP}) as
\begin{eqnarray}
\begin{array}{cl}
\min &  \c_B^{\T}\x_B + \c_N^{\T}\x_N,  \\
\mbox{\rm subject to} 
&  \A_B\x_B + \A_N\x_N=\b, 
\hspace{0.1in} \x_B \ge \0, \hspace{0.1in}  \x_N \ge \0.
\end{array}
\label{LP1}
\end{eqnarray}
Since $\A_B$ is non-singular, we can rewrite (\ref{LP1}) as
\begin{eqnarray}
\begin{array}{cl}
\min & \c_B^{\T}\A_B^{-1} \b 
+ (\c_N - \A_N^{\T} \A_B^{-\Tr} \c_B)^{\Tr} \x_N, \\
\mbox{\rm subject to} 
& \x_B = \A_B^{-1} \b - \A_B^{-1}\A_N\x_N, 
\hspace{0.1in} \x_B \ge \0, \hspace{0.1in}  \x_N \ge \0.
\end{array}
\label{LP2}
\end{eqnarray}
Let $\bar{\c}_N^{\Tr}=\c_N^{\T} -  \c_B^{\T} \A_B^{-1}\A_N$ be
the reduced cost, notice that the index of $\bar{\c}_N$ 
is the same as the index of $\A_N$. Denote the basic feasible solution
\begin{equation}
\bar{\b}=\A_B^{-1} \b =(\bar{b}_1,\ldots,\bar{b}_m) =\x_B.
\label{basicS}
\end{equation}
The first entering column $p$ of $\A_N$ is determined
by $\bar{c}_p = \min \{ \bar{\c}_N<0 \}$. 
Since $\a_t \in \A_N$ can be expressed as
$\a_t=\A_B{\y_t}$, therefore, we can write
\begin{equation}
{\y_t}=\A_B^{-1} \a_t=({y}_{t1},\ldots,{y}_{tm}).
\label{eColumn}
\end{equation}
The second entering column $q$ is determined by considering
all $\a_t$ corresponding to $c_t \in \bar{\c}_N<\0$ such that
\begin{eqnarray}
\bar{b}_j / y_{qj}=\max_{t \in \bar{\c}_N<0} 
\{ \min_{j \in \{ 1,\ldots,m \}} \bar{b}_j / y_{tj},   
\hspace{0.1in}  \mbox{subject to}  
\hspace{0.1in} y_{tj}>0 \}.
\label{outVar}
\end{eqnarray}
The leaving columns in $\A_B$ are obtained by solving
a two dimensional linear programming problem.
\begin{eqnarray}
\begin{array}{cl}
\min &  c_p x_p + c_q x_q,  \\
\mbox{\rm subject to} 
&  [\a_p ~~ \a_q] (x_p, x_q) =\bar{\b}, 
\hspace{0.1in} x_p \ge \0, \hspace{0.1in}  x_q \ge \0.
\end{array}
\label{LP2}
\end{eqnarray}
The indexes corresponding to vertex that achieves the 
optimum are the indexes of the leaving columns. 

%The following proposition shows this desired property.
%
%\begin{proposition}
%The cycling problem will not occur when the following
%pivot rules are adopted: (1) if the 
%longest step among all possible entering variables 
%is greater than zero, the longest step pivot is used;
%(2) if the longest step among all possible entering 
%variables is zero, the least index pivot rule is used.
%\label{firstP}
%\end{proposition}
%\begin{proof}
%First, we claim that there must exists a pivot rule that
%will bring the iterate out of degenerate basic feasible 
%solutions. This is known because both Bland's rule and 
%the lexicographic pivoting rule can achieve this. 
%Next, it is will known that the least pivot index
%rule will prevent cycling from happening. 
%Finally, a pivot that brings the iterate out of 
%degenerate basic feasible solutions must have a
%non-zero step size and the longest step size rule
%will achieve this goal. This proves the claim.
%\hfill \qed
%\end{proof}
%
%\begin{remark}
%The combined pivot rule described in the proposition
%may be more efficient than Bland's least index pivot
%rule because it may leave degenerate solution earlier
%than Bland's pivot rule. 
%\end{remark}

Because the double pivot algorithm achieves 
at least the same cost reduction as the longest step rule,
if the longest step rule results in a positive step size,
the double pivot rule will also achieve a positive
step size in the same iteration. This strategy
greatly increases the chance of bringing a non-degenerate 
variable (whose step length is greater than zero) into
the next basic feasible solution. We use Beale's problem
\cite{beale55} (Problem 2) to support this claim.
%\begin{problem}

\vspace{0.1in}
\noindent
{\bf Beale's problem} (Problem 2)
For this problem, $\c^{\T}=[-3/4,~150,~-1/50,~6,~0,~0,~0]$
\[
\A=\left[  \begin{array}{ccccccc}
1/4 & -60 & -1/25 & 9 & 1 & 0 & 0 \\
1/2 & -90 & -1/50 & 3 & 0 & 1 & 0 \\
0 & 0 & 1 & 0 & 0 & 0 & 1
\end{array}\right]
\hspace{0.1in} \b=\left[  \begin{array}{c}
0 \\ 0 \\1  \end{array}\right]
\] 
the initial base is $B^0=\{5,6,7\}$ and
$N^0=\{1,2,3,4\}$, $\bar{\c}_{N^0}^{\T}=[-3/4,~150,~-1/50,~6]$.
There are only two elements in $\bar{\c}_N^0<0$, the entering 
columns are $\{ 1,3 \}$. Therefore, we need to solve
\begin{eqnarray}
\begin{array}{cl}
\min &  -3/4 x_1 - 1/50 x_3,  \\
\mbox{\rm subject to} 
&  [\a_1 ~~ \a_3] (x_1, x_3) =\bar{\b}, 
\hspace{0.1in} x_1 \ge \0, \hspace{0.1in}  x_3 \ge \0.
\end{array}
\label{LP3}
\end{eqnarray}
Using graphic method, one can find the optimal solution of 
(\ref{LP3}) as $(x_1,~x_3)=(1/25,~1)$ with the last two rows 
composed of the vertex. Therefore, the leaving variables are 
$x_6$ and $x_7$. This gives $B^1=\{1,3,5\}$, $N^1=\{2,4,6,7 \}$,
and $\bar{\c}_{N^1}^{\T}=[15,~10.5,~1.5,~0.05]>\0$.
The optimal solution of this problem is found in one iteration.
While applying Dantzig's rule, after $6$ iterations, one will get 
$B^6=B^0$, cycling occurs.
%\end{problem}

For large $m$, an efficient algorithm was developed in 
\cite{ve18} to solve (\ref{LP2}).
More test shows that the double pivot algorithm finds the
optimal solutions for all cycling problems provided in this paper. 
%However, this does not show that this algorithm
%solves all cycling problems, our goal is to solve at least
%a majority of Netlib benchmark standard problems, if not all
%of these problems. 
The test result is summarized in the table below.
In this table, only basic feasible solutions and optimal basic
solutions are provided. All non-basic variables are zeros.

\tiny
\begin{longtable}{|c|c|c|c|c|}
\hline          
Problem & Initial point  & Iteration & Optimal solution & Objective function \\ 
\hline
1  & $[x_1, x_2, x_3]=[1, 0 ,0]$  & 2 & $[x_1,x_4,x_6]=[1,0,0]$ & 0 \\
\hline
2 &   $[x_5, x_6, x_7]=[0 ,0,1]$   & 1 &  $[x_1, x_3, x_5]=[0.04, 1, 0.03]$ & -0.05 \\
\hline
3 &  $[x_5, x_6, x_7]=[0 ,0,1]$   & 1 &  unbounded  & $-\infty$ \\
\hline
4 &   $[x_1, x_2, x_7]=[0 ,0,1]$   & 1 &  $[x_1, x_2, x_5]=[3, 2,1]$   & 1 \\
\hline
5 &  $[x_1, x_2, x_7]=[0 ,0,1]$    & 1 &  $[x_1, x_2, x_5]=[2.5, 1.5,0.5]$ & -0.5 \\
\hline
6 &  $[x_1, x_2, x_3]=[0 ,0,2]$    & 6 &  $[x_1, x_4, x_6]=[2, 2,2]$ & -2 \\
\hline
7 &  $[x_1, x_2, x_3,x_4]=[0,0,0,1]$    & 1 &  $[x_1, x_2, x_5,x_6]=[4, 1,4,1]$ & -2 \\
\hline
8 &  $[x_1, x_2]=[0,0]$    & 3 &  $[x_4, x_6]=[0,0]$ & 0 \\
\hline
9 &  $[x_5, x_6]=[0,0]$    & 1 &  unbounded &  $-\infty$  \\
\hline
10 &  $[x_1, x_7,x_8]=[2,0,2]$    & 1 & unbounded & $-\infty$ \\
\hline
11 &  $[x_5, x_6]=[0,0]$    & 4 &   $[x_2, x_5]=[0,0]$   &  0  \\
\hline
12 &  $[x_5, x_6]=[0,0]$    & 1 & unbounded & $-\infty$ \\
\hline
13 &  $[x_5, x_6]=[0,0]$    & 4 &   $[x_2, x_5]=[0,0]$   &  0  \\
\hline
14 &  $[x_5, x_6, x_7,x_8]=[0,0,5,10]$    & 1 &  
$[x_1, x_3, x_5,x_7]=[0.0474, 3.3175,0.2844,1.6351]$ & -1.8246 \\
\hline
15 &  $[x_7, x_8]=[0,0]$    & 1 &   unbounded  &   $-\infty$  \\
\hline
16 &  $[x_9, x_{10}]=[0,0]$    & 2 &   unbounded  &   $-\infty$  \\
\hline
17 &  $[x_5, x_6]=[0,0]$    & 1 & unbounded & $-\infty$ \\
\hline
18 &  $[x_1, x_2, x_3]=[0,0,0]$    & 7 &  $[x_2,x_3,x_9]=[0,0,0]$ & 0 \\
\hline
19 &  $[x_1, x_2]=[0,0]$    & 4 &  $[x_1,x_5]=[0,0]$ & 0 \\
\hline
20 &  $[x_5, x_6, x_7,x_8]=[0,0,1,1]$    & 2 &  
$[x_1, x_2, x_4,x_6]=[1, 1,3,8]$ &  -3.2500 \\
\hline
21 &  $[x_5, x_6]=[0,0]$    & 1 & unbounded & $-\infty$ \\
\hline
22 &  $[x_5, x_6]=[0,0]$    & 1 & unbounded & $-\infty$ \\
\hline
23 &  $[x_7, x_8]=[0,0]$    & 2 & unbounded & $-\infty$ \\
\hline
24 &  $[x_9, x_{10}]=[0,0]$    & 1 & unbounded & $-\infty$ \\
\hline
25 &  $[x_9, x_{10}]=[0,0]$    & 1 & unbounded & $-\infty$ \\
\hline
26 &  $[x_8, x_9, x_{10}]=[0,0,0]$    & 6 & unbounded & $-\infty$ \\
\hline
27 &  $[x_5, x_6, x_7]=[0,0,1]$    & 3 & unbounded & $-\infty$ \\
\hline
28 &  $[x_5, x_6, x_7]=[0,0,1]$    & 5 & $[x_1, x_3, x_5]=[0.5,0.5,1]$ & $-0.5$ \\
\hline
29 &  $[x_5, x_6, x_7]=[0,0,0]$    & 3 & unbounded & $-\infty$ \\
\hline
30 &  $[x_7, x_8, x_9]=[0,0,0]$    & 3 & unbounded & $-\infty$ \\
\hline
\caption{Numerical test for double pivot method}
\label{tableIteration}
\end{longtable}
\normalsize

\section{Conclusions}
In this paper, we collected a set of cycling problems in linear
programming. This set of problems is used to test a double 
pivot simplex algorithm. The test result indeed shows that 
the double pivot algorithm prevents the cycling from happening 
for all the problems presented in this paper. Our ultimate goal 
in the next project is to show that this double pivot algorithm,
together with some other simple strategies, will be able
to solve a majority set of very challenging Netlib benchmark
problems, if not all them. These problems can be downloaded 
from Matlab file exchange cite https://www.mathworks.com/matlabcentral/fileexchange/91985-a-collection-of-cycling-problems-in-linear-programming.

\section{Acknowledgment}
This author thanks Professor Zörnig for kindly providing his
paper that is very helpful in the preparation of this paper.

\vspace{0.4in}
Declarations of interest: none.

%\section{Conflict of interest}
%
%On behalf of all authors, the corresponding author states that 
%there is no conflict of interest.

\end{document}